\input amstex.tex
\documentstyle{amsppt}
\magnification\magstep1
\NoBlackBoxes
\define\tr{\text{tr}}
\NoRunningHeads
\topmatter
\title THE CHERN CHARACTER OF CERTAIN INFINITE RANK BUNDLES ARISING IN GAUGE THEORY \endtitle
\author Jouko Mickelsson \endauthor
\affil Department of Mathematics and Statistics, University of Helsinki, and Royal Institute of Technology, Stockholm \endaffil
\endtopmatter 
\hoffset=1cm
\voffset=1cm
\advance\vsize -1cm

\document

ABSTRACT:  A cocycle  $\Omega: P \times G \to H$ taking values in a Lie group $H$ for a free right action of $G$ on $P$ defines a principal
bundle $Q$ with the structure group $H$ over $X=P/G.$ 
The Chern character of a  vector bundle associated to $Q$ defines then characteristic classes on $X.$ This observation becomes useful in the
case of infinite dimensional groups. It typically happens that a representation of $G$ is not given by linear operators which differ from the identify
by a trace-class operator. For this reason the Chern character of a vector bundle associated to the principal fibration $P \to P/G$ is ill-defined. 
But it may happen that the Lie algebra representations of the group $H$ are given in terms of trace-class
operators and therefore the Chern character is well-defined; this observation is useful especially if the map $g\mapsto \Omega(p;g)$ is a homotopy
equivalence on the image for any $p\in P.$ We apply this method to the case $P= \Cal A,$ the space of gauge connections
in a finite-dimensional vector bundle, and $G= \Cal G$ is the group of (based) gauge transformations. The method for constructing
the appropriate cocycle  $\Omega$ comes from ideas in quantum field theory, used to define the renormalized gauge currents in a Fock space.

\vskip 0.3in

1. INTRODUCTION
 
\vskip 0.2in

In gauge theory one is led to study infinite rank vector bundles where the structure group is the group $U^p,$  resp. $U_p,$ (or its subgroup), consisting of
unitary operators $g$ in a complex Hilbert space $H$ such that $g-1$ is in $L_p$ in the even case, resp  $[\epsilon, g]\in L_p$ in the odd case, [MR]. 
Here $\epsilon$ defines a polarization in $H,$ and even/odd refers to the parity of the dimension of the underlying compact physical space $M.$
$L_p$ is the Schatten ideal of bounded operators $T$ such that $|T|^p$ is trace-class, $p\geq 1.$ The index $p$ depends on the dimensionality 
of $M,$ $p > n=\text{dim}\, M.$  The reason for the choice of $p$ is related to the fact that the Green's function $D^{-1}$ of a Dirac operator $D$
is in $L_p$ for $p> n.$ 

Starting from the even case, one can pose the more general question how to define the Chern character for $U^p$ bundles over some space $X.$ Typically, $X$ is the
moduli space of gauge connections. The problem in the infinite rank setting is that the powers $\Cal F^k$ of the curvature are not always trace-class 
operators. That would be the case in the even case if $k  \geq p.$  \newpage

On the other hand, all the groups $U^p$ are homotopy equivalent to $U^1,$ 
due to a theorem by Richard Palais [P].  Concretely, the $U^p$ bundles arise from principal gauge group bundles through the embedding
of the group of gauge transformations $\Cal G$ to $U^p.$  The group $U^p$ can the be replaced by $U^1$ using a cocycle $\Omega: P \times\Cal G
\to U^1.$ For each point $A\in P$ the map $\Cal G \to U^1$ defined by $g\mapsto \Omega(A;g)$ is homotopic to the embedding $\Cal G \subset U^p,$
in the topology of the larger group $U^p.$  This method has been previously used in the quantization of chiral fermions in background gauge fields 
and in computation of hamiltonian anomalies (Schwinger terms), [M]. 

In order to define the Chern character one needs a connection. Normally, when a vector bundle is associated to a principal bundle $P$ over $X$ the 
connection form acting on sections of the vector bundle is obtained from the principal bundle connection by the action of the Lie algebra representation $\rho$
in the model fiber of the vector bundle. However, in the case when the associated vector bundle is defined by a cocycle $\Omega$ there is a complication.
There are two candidates for a connection. The first comes from
the canonical 1-form defined  by the derivative $\omega$ of $\Omega$ with respect to the second argument,   and computing $\omega(p; \theta)$ with $\theta$ a 
connection 1-form for the bundle $P\to X,$ is a natural candidate to replace $\rho,$ but this does not quite define
a connection.   It does not  transform properly under the right translations by $\Cal G.$ The second candidate is determined from the 'renormalization process'
defining the cocycle $\Omega.$ That does define a covariant differentiation in the associated vector bundle. However, the connection form takes values
in a larger Lie algebra than $\bold u^1 = Lie(U^1).$ It turns out that we can still use this latter 1-form to compute the characteristic classes. The reason
for this is that the gauge variations of the corresponding Chern-Simons forms are given by traces of trace class operator valued forms, due to the fact
that the transition functions take values in $\bold u^1.$  Pairing Chern classes with singular homology classes and using Stokes'  theorem the
integral of a Chern class against a singular cycle becomes a sum of integrals of the gauge variations of Chern-Simons form over the faces of
the individual singular simplices  and thus are well-defined.

\vskip 0.3in

2. CHERN CHARACTER FOR $U^1$ BUNDLES AND FAMILIES OF DIRAC OPERATORS

\vskip 0.2in

Let $M$ be a  compact  oriented connected spin manifold and $G$ a compact Lie group and $\frak g$ its Lie algebra. 
 Denote by $\Cal A$ the space of smooth $\frak g$ valued 1-forms on $M,$ 
the gauge potentials in a trivial $G$ bundle over $M.$  In the following we could generalize the discussion to nontrivial $G$ bundles, but I 
want to avoid unnecessary technical complications which actually are irrelevant for the main ideas below. Denote by $\Cal G$ the group of 
gauge transformations (bundle automorphisms which fix the base), which in the case of a trivial bundle can be identified as the group
of smooth maps from $M$ to $G.$  

Let us first consider the case when the dimension of $M$ is odd.
In quantum field theory one wants to associate a bundle of fermionic Fock spaces to $\Cal A,$ such that the quantized Dirac operators $\hat D_A$ with
$A\in\Cal A,$ transform covariantly under the group of gauge transformations. It turns out that in the case of chiral Dirac operators (acting on massless Weyl fermions)
there is generically 
an obstruction to this, in the form of a \it commutator anomaly \rm (Schwinger terms) which means that instead of $\Cal G$ one has to construct certain
group extension $\hat\Cal  G$ which acts on sections of the Fock bundle, [M].  There is a technical difficulty in this construction: The group
$\hat\Cal G$ has an unitary representation in the Fock space only when dim$\,M =1.$ 

When the dimension is higher than one the representations
have to be understood in a generalized sense. Instead of a group homomorphism $\hat\Cal G \to U(\Cal H),$ where $\Cal H$ is the fermionic Fock space, 
one has to construct first a cocycle $\Omega: \Cal A \times \Cal G  \to U(H),$ where $H$ is the fermionic 1-particle space, such that the values of $\Omega$
are in the restricted unitary group $U_{res}(H) \subset U(H).$ The restricted unitary group consists of unitaries such that the off-diagonal blocks with 
respect to a fixed polarization of $H$ are Hilbert-Schmidt operators.  In that case the operators $\Omega(A;g)$ can be represented by unitaries 
$\hat\Omega(A;g) \in U(\Cal H)$ which then define a group extension of $\Cal G.$   

In the present article I want to consider the problems arising in the case of an even dimensional manifold $M$ from the point of view of families index theory.
In short, the problem is the following: Consider $\Cal A$ as a principal bundle over the moduli space $\Cal A/\Cal G$ with fiber $\Cal G.$ As $\Cal G$ 
I take now the group of smooth \it based maps \rm from $M$ to $G.$ Based means that we fix a point $x_0\in M$ and consider maps $g$ such that
$g(x_0)=1,$ the neutral element in $G.$ This guarantees  that $X=\Cal A/\Cal G$ has no singularities, it is a Frechet manifold.  Given 
a finite dimensional representation $\rho:G \to \text{Aut}(V)$ we can construct a Hilbert bundle $\Cal E$ over $\Cal A/ \Cal G$  through the natural action of $\Cal G$ 
on sections of $S(M) \otimes E$ where $E$ is a (trivial) complex vector bundle over $M$ with fiber $V$ and $S(M)$ is the spin bundle over $M.$ 

Any Hilbert bundle with structure group $U(H)$ (in the  norm topology) is topologically trivial by Kuiper's theorem, [K].  However, in our case there is additional geometrical structure. The bundle has $\Cal G$ as its 
structure group and we can define a connection on $\Cal E$ taking values in the Lie algebra of $\Cal G.$ The problem next is to define the Chern character
of the bundle in the usual way as $\tr\, e^{\Cal F/2\pi i},$ where $\Cal F$ is the curvature of the connection. As it stands, the Chern character is ill-defined since 
the trace of the curvature diverges. One way to avoid the divergence is to replace the curvature by the curvature of an appropriate Quillen superconnection, [Q].
This method was worked through in detail by J. - M. Bismut, [B], for the derivation of the local Atiyah-Singer index theorem for families of Dirac
operators, [AS1].  In this article I want to propose a different solution, very close to the idea of renormalization
in [M] in the odd dimensional case. 

Let us first assume that the Dirac index of $D_A$ is zero. Then we can fix an invertible  hermitean operator $D_0$ such that $D_A -D_0$ is bounded for any
$A\in\Cal A.$ Let $\epsilon = D_0/|D_0|.$  Perturbing $D_A$ by a suitable finite rank operator $q$ one obtains an invertible operator $D'_A = D_A +q$ which 
anticommutes with the grading operator $\Gamma$ and so also $F_A = D'_A/|D'_A|$ anticommutes with $\Gamma.$ 

Using the method in [LM] one can construct a family $T_A$ of unitary operators commuting with $\Gamma$ such that ${T_A}^{-1} F_A T_A - \epsilon$ is a trace-class
operator. The operators $T_A$ are not uniquely defined, if $R_A \in U_1$ (with respect to the polarization defined by the sign $\epsilon$) then also $R_A T_A$ 
satisfies the condition. In order to be self contained, we give the (somewhat simplified) proof adapted to the present discussion.

\redefine\e{\epsilon}

\proclaim{Theorem 1 } For all interactions $A$ in $D_A = D_0 + A$ with $A$ a pseudodifferential operator of order zero, there is an
unitary pseudodifferential operator $T_A,$  commuting with the spinor grading operator $\Gamma,$ such that
the transformed  interaction $A'$
satisfies  $[\e,A'] \in L_1$ and $T_A =1$ for $A=0.$  The symbol of $T_A$ depends smoothly on the symbol of $A.$ \endproclaim

\demo{Proof of Theorem}
We first prove the following
key lemma providing the recipe for constructing $T_A$:

\proclaim{Lemma 1} Let $A$ be a pseudodifferential operator such that
$[\e,A]$ is of order $k \leq 0.$  Then $A',$ defined by $D_0 + A' = T_A^{-1} (D_0 +A) T_A $ 
 with
the unitary operator $$ T_A =  e^{\alpha} ,\quad \alpha=
-\frac18\left(|D_0|^{-1}[\e,A]+[\e,A]|D_0|^{-1}\right)  $$
has the property  $[\e,A']$ is of order $k-1$ 
\endproclaim

\demo{Proof of Lemma}  Denote the space of pseudodifferential operators of order $k$ by $O_k;$ this is in fact an algebra if $k\leq 0.$
We write $A' =  A'_1 +
A'_2$ where
$$A'_1 =  A + [D_0,\alpha] $$
is the leading term
in the asymptotic  expansion, and
$$A'_2 =   T_A^{-1}[D_0,T_A-\alpha-1] +  T_A^{-1}[A,T_A-1] $$
is the rest.

Since obviously $\alpha$ is  of order $k-1$, and the fact that 
$A'_2$ is of order $k-1$ trivially follows from $\alpha D_0, D_0\alpha \in O_{k}$
and 
$$ T_A-1 -\alpha = \frac12 \alpha^2 + \frac{1}{6} \alpha^3 +\dots \in O_{2k-2} \subset  O_{k-2}.  $$

The nontrivial part thus is to show that  $[\e,A'_1]\in O_{k-1}$. This
can be seen by the following calculation,
$$
\align
[\e,A_1'] =  [\e,A]  -\frac{1}{8}\left[\e,\left[D_0,|D_0|^{-1}[\e,A] +
[\e,A]|D_0|^{-1}\right] \right] \\ =
\frac18\left(8[\e,A] -\left[\e, \e[\e,A] - |D_0|^{-1}[\e,A]|D_0|\e
+ \e |D_0|[\e,A]|D_0|^{-1} - [\e,A]\e \right] \right) \\=
\frac18\left(8[\e,A] - 4[\e,A] - 2|D_0|[\e,A]|D_0|^{-1} -
2|D_0|^{-1}[\e,A]|D_0| \right) \\=
\frac14\left(|D_0|^{-1}\left[|D_0|,[\e,A]\right] - \left[|D_0|,[\e,A]\right]
|D_0|^{-1} \right)
\endalign
$$
where we used $[\e,A] =  -\e [\e,A]\e$ and
$\e D_0^{\pm 1} =  D_0^{\pm 1}\e =  |D_0^{\pm 1}|$.
Thus
$$[\e,A_1']  =
\frac{1}{4}\left[|D_0|^{-1},\left[|D_0|,[\e,A]\right]\right]
$$ which is in $O_{k-1}$ since $[ |D_0|, O_{p}] \subset O_{p}.$  If we replace $D_0^{-1}$ by $D_0(D_0^2+\lambda)^{-1}$
a similar calculation leads to the same conclusion,
$$
[\e,A_1'] =  \frac14\left[
|D_0|(D_0^2+\lambda)^{-1},\left[|D_0|,[\e,A]\right] \right] +
\frac12\lambda\left((D_0^2+\lambda)^{-1}[\e,A]+[\e,A](D_0^2+\lambda)^{-1}\right).
$$
This proves our Lemma.
\enddemo

We can apply this method successively: Starting from some interaction $A_0=
A$ such that $[\e,A]\in O_0$ we get a new interaction $A_1=A'$ using the
conjugation $T_A,$ with $[\e,A_1]\in O_{-1}.$ We can then insert $A_1$ as
an argument to $T$ and obtain an unitary operator $T_{A_1}.$ This
defines again a new interaction $A_2=A_1'$ such that $[\e,A_2]\in O_{-2}.$
Continuing this way we obtain, after $p$ steps, an unitary operator
$T^p_A = T_{A_{p-1}}\dots T_{A_0}$ such that the interaction $A_p$ in
 $T^p_A (D_0 +A)(T^p_A)^{-1}=D_0 + A_p$  has the property
$[\epsilon, A_p] \in L_1.$ The number $p$ of iterations  depends on the dimension of $M,$ one must take
$p > \text{dim}\, M.$       $ \square$
\enddemo

\bf Remark \rm  Actually, in quantum field theory applications (and also in Section 3 below) it is sufficient that $[\epsilon, A']$ is Hilbert-Schmidt;
in practice, this means that one needs to perform one less iteration in the algorithm.

Since $T_A$ commutes with $\Gamma,$ the operator ${T_A}^{-1} F_A T_A$ anticommutes with $\Gamma$ and we can write
$$ {T_A}^{-1} F_A T_A = \left( \matrix 0 & g\\ g^* & 0 \endmatrix\right)$$
where the blocks are defined with respect to the grading $\Gamma.$ Here $g: S_- \to S_+$ is a unitary operator from the left chirality spinors to right 
chirality spinors.  Likewise,  
$$\epsilon = \left( \matrix 0 & g_0 \\ g_0^* & 0\endmatrix\right).$$
with $g - g_0$ a trace-class operator. Thus we can write $g_0^{-1} g = h_A$ and $h_A-1$ is  trace-class. If $T'_A$ is another unitary
operator which satisfies the condition $T_A^{-1} F_A T_A - \epsilon \in L_1$ then we can write (by the ideal property of $L_1$) $T'_A =  T_A S_A$
where $S_A \in U_1.$ 

Since $T_A$ and the gauge transformations $h\in \Cal G$ commute with $\Gamma,$ we can write
$${T_A}^{-1} h T_{A^h} = \Omega(A;h) = \left(\matrix \Omega_+(A;h) & 0\\ 0 & \Omega_-(A;h) \endmatrix\right)\tag2.1$$
where $\Omega$ satisfies the cocycle relation
$$\Omega(A;h_1)\Omega(A^{h_1};h_2) = \Omega(A; h_1h_2).\tag2.2$$
The change $T_A\mapsto T'_A = T_A S_A$ acts on the cocycle as
$$ \Omega'(A; h) = {S_A}^{-1} \Omega(A;h) S_{A^h} .\tag 2.3$$
 
\proclaim{Lemma 2} The cocycle $\Omega= (\Omega_+, \Omega_-)$ satisfies
$$ \Omega_+ - g_0 \Omega_- g_0^{-1}  \in L_1.$$
\endproclaim
\demo{Proof}   By ${T_{A^h}}^{-1} F_{A^h} T_{A^h} - \epsilon \in L_1$ and $F_{A^h} = h^{-1} F_A h + \text{finite rank operator},$
we obtain by 
$$hT_{A^h} \epsilon {T_{A^h}}^{-1} h^{-1} -F_A \in L_1$$
and subtracting from $T_A \epsilon {T_A}^{-1}  - F_A \in L_1$ the relation
$$h T_{A^h} \epsilon {T_{A^h}}^{-1} h^{-1} -T_A \epsilon {T_A}^{-1} \in L_1.$$
Finally, multiplying this relation from the left by ${T_A}^{-1}$ and from the right by $hT_{A^h}$ gives   
$$[\epsilon, {T_A}^{-1} h T_{A^h}] = [\epsilon, \Omega(A; h)] \in L_1,$$
which  can be written as a pair of relations
$$g_0 \Omega_-  - \Omega_+ g_0 \in L_1 \text{ and } g_0^{-1} \Omega_+ - \Omega_- g_0^{-1} \in L_1,$$ 
which actually are equivalent. The claim follows immediately by multiplying the second relation by $g_0$ from the left  and
from the ideal property of $L_1.  \square$ \enddemo

\proclaim{Corollary 1} The conditional supertrace $\text{str}_c (\Omega) = \frac12 \text{tr}(\Gamma \Omega +\epsilon \Gamma\Omega \epsilon)$ is absolutely converging.
\endproclaim
If $\Omega$ is trace-class then the conditional supertrace above is just the supertrace $\text{tr}\, \Gamma\Omega.$ 

The set of unitary matrices $\Omega= \left( \matrix \Omega_+  & 0 \\ 0 & \Omega_-\endmatrix\right)$ satisfying the property in Lemma 2 is an infinite dimensional
Lie group $UU^1(H)$ which is homotopy equivalent to $U^1.$ Namely, its elements are parametrized by the unitaries $\Omega_-$ acting in the Hilbert space
of left handed spinor fields and the unitary matrix $z = (g_0 \Omega_-)^{-1} \Omega_+ g_0 .$ By repeated use of Lemma 2 the operator $z$ acting on left-handed 
spinors differs from the identity by a trace-class operator.

\bf Remark \rm Since the space of gauge connections is contractible, $(A,s) \to sA$ with $0\leq s \leq 1$ is a contraction, the cocycle $\Omega(A;h)$ is 
contractible to the tautological cocycle $(A,h) \mapsto h = \Omega(sA;h)$ for $s=0.$ However, this homotopy is defined only in the framework of
$U(H)$ valued cocycles, not in cohomology with values in $UU^1(H).$ 

The above construction can be extended to the case when ind$D_A \neq 0.$ Now we fix $\epsilon$ such that the parameter $g_0: S_- \to S_+$ has Fredholm
index equal to the Dirac index of $D_A.$ In addition, we can require that $g_0$ is an isomorphism from the orthogonal complement of its kernel to the orthogonal
complement of the
cokernel. Likewise, one can choose $\epsilon_A$ such that it is equal to the sign of $D_A$ in the orthogonal complement of ker$D_A$ and is equal to zero
on ker$D_A.$ Then $\epsilon_A$ anticommutes with $\Gamma.$  The construction of $T_A$ is essentially the same as in the index zero case: Instead of
$D_0^{-1}$ in the formula for $T_A$ one chooses an operator $R_0$ such that $R_0 D_0$ is equal to the identity on the orthogonal complement of ker$D_0.$
This does not fix $R_0$ uniquely, only up to a finite rank operator, but that is sufficient for the construction of the cocycle $\Omega$ in the Lemma. 

The cocycle $\Omega$ defines a principal bundle $P$ over $X$ with structure group $UU^1(H).$ The total space of $P$ is 
$$P = \Cal A \times_{\Omega} UU^1(H,)$$
the set of equivalence classes of pairs $(A, g) \in \Cal A \times UU^1(H)$ where the equivalence relation is $(A,g) \equiv (A^{h^{-1}}, \Omega(A; h)g)$
where $h\in \Cal G.$ The right action of $UU^1(H)$ on $P$ is given by the right multiplication on the second component. 

Given a system $\psi_{\alpha}: U_{\alpha} \to \Cal A$ of local trivializations of the bundle $\Cal A \to X$ one has local trivializations 
$\tilde\psi_{\alpha} : U_{\alpha} \to P,$ $\tilde\psi_{\alpha}(x) = (\psi_{\alpha}(x),1).$ The transition functions are  
$\tilde g_{\alpha\beta}(x)= \Omega(\psi_{\alpha}(x), g_{\alpha\beta}(x))$ and they satisfy the usual cocycle relation
$$\tilde g_{\alpha\beta} \tilde g_{\beta\gamma} = \tilde g_{\alpha\gamma}$$
as a consequence of the twisted cocycle relation for $\Omega.$

Following [AS2], a connection in the principal bundle $\Cal A \to \Cal A/ \Cal G$ is constructed as follows. For $A\in \Cal A$ the horizontal subspace
$H_A$ in the tangent space of $\Cal A$ consists of all $B \in \Cal A$ such that $\sum \partial_i B^i + [A_i, B^i]=0.$ Here $B^i = \sum g^{ij} A_j$ where 
$g^{ij}$ is the inverse of the metric tensor $g_{ij}$ on $M.$  If $B\in H_A \cap V_A$ then on one hand $B_i =D_i^A X$ for some $X\in C^{\infty}(M, \frak g)$
but on the other hand $D^A_i B^i = \Delta_A X =0.$  The harmonic functions $X$ satisfy
$$ 0 = \int_M <X, \Delta_A X> dvol_M = \int_M < D^A_i X, D^{A, i} X >,\tag2.4$$
where $<\cdot,\cdot>$ is an invariant inner product on $\frak g.$ Thus $D^A_i X=0$ for each coordinate index $i$. On a connected manifold $M$ a covariantly
constant function $X$ is completely fixed by its value at a base point $x_0.$ Thus for based infinitesimal gauge transformations we must have $X=0$ everywhere
and so $V_A \cap H_A =\{0\}.$ 

The connection form $\theta$ corresponding to the horizontal distribution $H_A$ can be given fairly explicitly, [AS2]. Since the only harmonic function vanishing
at the base point is zero,  the form
$$\theta_A(B) = \Delta_A^{-1} (D^A_i B^i)\tag2.5$$
is uniquely defined, with values in the Lie algebra $Map_0(M,\frak g).$ It vanishes on horizontal forms $B$ and is tautological on the vertical tangent vectors
$B_i = D^A_i X.$ 

Given a principal bundle $\pi:P \to M$ with structure group $G$ and a cocycle $\Omega : P \times G \to K$ where $K$ is a group of invertible bounded operators
in a Hilbert space $V$, one can construct
an associated vector bundle $E= P \times_{\Omega} V.$ The elements of $E$ are equivalence classes of pairs $(p,v)\in P \times V$ with the equivalence
relation defined by $(p,v) \sim (pg^{-1}, \Omega(p;g)v).$  As in the standard case, when the cocycle $\Omega$ is simply a homomorphism $\rho(g) =
\Omega(p;g),$ a connection in the principal bundle $P$ should define  a connection in the associated vector bundle. The construction is slightly more
complicated 
in the more general situation. A  section $\xi$ of $E$ can be identified as a
map $\xi: P \to V$ such that $\xi(pg) = \Omega(p;g)^{-1} \xi(p)$ for all $(p,g)\in P\times G.$ Fixing a local section $\psi$ of $P,$ a section of $E$ is given
by a local vector valued function $f = \xi\circ \psi.$  Denote by $\omega$ the derivative of $\Omega(p;g)$ with respect to the second argument at $g=1.$ 

Returning to the case when $P=\Cal A$ is the space of gauge connections and $\Omega$ is the cocycle (2.1) we have
\proclaim{Lemma 3} A change $T'_A = S_A T_A$ induces a gauge transformation
$$ \omega'(A; \theta(v)) = S_A \omega(A;\theta(v)) {S_A}^{-1} - (\Cal L_{\theta(v)} S_A) {S_A}^{-1} $$
on the  form $\omega.$ \endproclaim
\demo{Proof}  Differentiating (2.3) with respect to parameter $t$ at $t=0$ for $h=e^{t\theta(v)}$ we obtain
$$\align \omega'(;\theta(v)) &  = \frac{d}{dt}|_{t=0} \Omega'(A; e^{t\theta(v)}) = \frac{d}{dt}|_{t=0} S_A\Omega(A; e^{t\theta(v)}) {S_{A^h}}^{-1} \\
& = S_A \omega(A;\theta(v)){S_A}^{-1}  - (\Cal L_{\theta(v)} S_A) {S_A}^{-1}\endalign $$
where $v$ is a tangent vector on $\Cal A$ at the point $A$ and $\Cal L_{\theta(v)}$ is the Lie derivative in the direction of the infinitesimal gauge transformation $\theta(v).$ 
 $\square$ \enddemo

The natural form $\omega(p;\theta)$ taking values in the Lie algebra of $K$  is not quite a connection form as the following proposition shows:

\proclaim{Proposition 1} Let $\omega_{\alpha}$ be the pull-back of the form $\omega(p; \theta)$ with respect to the local section $\psi_{\alpha}$ of the
bundle $P\to M.$  Then $$\omega_{\beta} = \Omega(\psi_{\alpha} ; g_{\alpha\beta})^{-1} \omega_{\alpha} \Omega(\psi_{\alpha}; g_{\alpha\beta}) 
+ \Omega(\psi_{\alpha}; g_{\alpha\beta})^{-1} D^{vert,1} \Omega(\psi_{\alpha}; g_{\alpha\beta})\tag2.6$$
where $D^{vert,1}$ means the differentiation  taken in the vertical directions in the first argument,
$$D^{vert,1}_v \Omega(\psi; g)=  \frac{d}{dt} \Omega( \psi\cdot e^{t\theta_{\alpha}(v)} ; g)\vert_{t=0} + D^{(2)}_v \Omega(\psi_{\alpha}; g),$$
where $D^{(2)}_v$ is the differentiation applied to the second argument, in the direction of the vector $v,$  and $\theta_{\alpha}=
\psi_{\alpha}^* \theta$ is the pull-back to $U_{\alpha}.$ \endproclaim

\demo{Proof} Using the cocycle identity $\Omega(p;g_0) \Omega(p\cdot g_0; {g_0}^{-1} g) = \Omega(p;g)$ we obtain first
$$D^{(2)}_v \Omega(p;g) = D^{(2)}_v \Omega(p;g_0 {g_0}^{-1} g)= \Omega(p;g_0) D^{(2)}_v \Omega(p\cdot g_0; {g_0}^{-1} g)$$
where $D^{(2)}_v$ (resp. $D^{(1)}_v$) denotes the derivative with respect to the second (resp. first) argument in the direction of the vector
$v$ on $M;$ here $g_0$ is a constant and $g=g(z)$ is a group valued function on $M.$

Taking now $g=g_{\alpha\beta}(z)$ and $p=\psi_{\alpha}(z),$ and after taking the derivatives set $g_0 = g_{\alpha\beta}(z),$ we obtain
$$ \omega(\psi_{\alpha}(z)\cdot g_{\alpha\beta}(z) ; {g_{\alpha\beta}}^{-1} D_v g_{\alpha\beta})
= \Omega(\psi_{\alpha}(z) ; {g_{\alpha\beta}}(z) )^{-1} D^{(2)}_v \Omega(\psi_{\alpha}(z) ; g_{\alpha\beta}(z)).\tag2.7$$
Next, using again the cocycle property in the second equality below, we obtain
$$ \align & \omega(\psi_{\beta}(z); {g_{\alpha\beta}(z)}^{-1} \theta_{\alpha}(v) g_{\alpha\beta}(z))  \tag2.8\\
 & = \frac{d}{dt}|_{t=0} \Omega(\psi_{\alpha}(z) \cdot g_{\alpha\beta}(z) ; {g_{\alpha\beta}(z)}^{-1} e^{t\theta_{\alpha}(v)} g_{\alpha\beta}(z)) \\
& = \frac{d}{dt}|_{t=0} \Omega(\psi_{\alpha}(z) \cdot g_{\alpha\beta}(z) ; {g_{\alpha\beta}(z)}^{-1} ) \Omega(\psi_{\alpha}(z) ; e^{t\theta_{\alpha}(v)} g_{\alpha\beta}(z)) \\
& =\frac{d}{dt}|_{t=0} \Omega(\psi_{\alpha}(z); {g_{\alpha\beta}(z)} )^{-1}\Omega(\psi_{\alpha}(z) ; e^{t\theta_{\alpha}(v)} g_{\alpha\beta}(z)) \\
& = \frac{d}{dt}|_{t=0} \Omega(\psi_{\alpha}(z) ; g_{\alpha\beta}(z))^{-1} \Omega(\psi_{\alpha}(z); e^{t\theta_{\alpha}(v)}) \Omega(\psi_{\alpha}\cdot e^{t\theta_{\alpha}(v)} ; g_{\alpha\beta}(z))\\
& = \Omega(\psi_{\alpha} ; g_{\alpha\beta})^{-1} \omega(\psi_{\alpha}; \theta_{\alpha}(v)) \Omega(\psi_{\alpha}  ; g_{\alpha\beta})
+ \Omega(\psi_{\alpha}; g_{\alpha\beta})^{-1} D^{(1,vert)}_v \Omega(\psi_{\alpha}; g_{\alpha\beta}). 
 \endalign  $$
 
 Adding equation (2.8) to (2.7) we arrive at (2.6) 
$\square$ 

Thus the defect in $\omega(p; \theta),$ for being a connection form, is the missing differentiation in the horizontal direction in the first argument on the
right-hand-side in (2.6). A true connection form $\theta_Q$  on  the principal $K$ bundle $Q= P\times_{\Omega} K$ comes from a $\bold k$ valued 1-form $\theta_P$ on $P$
as $\theta_Q = k^{-1} \theta_P k + k^{-1} dk$ where $k\in K$ if and only if $\theta_P$ satisfies the following conditions:
\roster
\item  $\theta_P^g = \Omega(\cdot; g)^{-1} \theta_P \Omega(\cdot; g)  - \Omega(\cdot; g)^{-1} d_P \Omega(\cdot;g)$
\item  $\theta_P(p; \hat X) = \omega(p; X)$
\endroster
where $\hat X$ denotes the vertical vector field on $P$ generated by the element $X$ in the Lie algebra of $G.$  These two conditions guarantee that
1) the form $\theta_Q$ on $P\times K$ indeed descends to $P\times_{\Omega} K$, 2) it is tautological  in vertical $K$ directions, and 3) is equivariant with
respect to to the combined right and adjoint action of  $K.$  In the more standard situation, when $\Omega$
does not depend on the argument $p\in P,$ the conditions above are the usual properties of a connection form. 

A true connection exists always in the
finite dimensional case. In an infinite dimensional case,  for a non-Hilbert manifold, one has to prove the existence of a partition of unity subordinate to an 
open trivializing  cover of the base manifold. In the present case when $P= \Cal A$ and $X=\Cal A/\Cal G$ such a partition of unity exists as proven in
[CMM]; the proof is based  a theorem of John Milnor [Mi].  So given an open cover $\{U_{\alpha}\}$ of $X$ with local trivializations $\psi_{\alpha}$ and 
transition functions $\tilde g_{\alpha\beta}$ the standard construction applies, the local connection forms can be defined as
$$ -\sum_{\beta} \rho_{\beta}  d \tilde g_{\alpha\beta} \tilde g_{\alpha\beta}^{-1}. $$
However, this definition is not very user friendly for the calculation of characteristic classes because of the necessity to pick up the partition of unity.

In the case when $P=\Cal A$  and $\Omega(A; g) = T_A^{-1} g T_{A^g}$ we can construct a connection as a covariant derivative acting on sections of the
associated vector bundle through the global 1-form $\tilde\theta= T_A^{-1} \theta T_A + T_A^{-1} d T_A.$ Although this form  does  not take values in the Lie 
algebra of the group $UU^1$ the covariant differentiation of a section $f$ of the vector bundle can nevertheless be defined as
$$\nabla f_{\alpha}  =  d f_{\alpha} + \tilde\theta_{\alpha} f.$$
We observe from the formula  $\Omega(A;g)= T_A^{-1} g T_{A^g}$ that actually $\tilde\theta$ is equal to $\omega(A; \theta)$ \it along the vertical 
directions \rm of the fibration $\Cal A \to \Cal A/ \Cal G=X.$ Along horizontal directions the latter form is zero whereas the former is nonvanishing. 
It might sound strange that the true connection is nonvanishing along horizontal directions but one should bear in mind that the horizontal 
directions in the fibration $\Cal A / \Cal G$ are not the same as the horizontal directions in the associated $UU^1$ principal bundle over $\Cal A / \Cal G.$ 
This fact is due to the nonvanishing of the derivative of $\Omega$ with respect to the first argument in the horizontal directions of the fibration
$\Cal A \to X.$

Although the curvature $\Cal F$ constructed from the connection form $\tilde\theta$ does not take values in the space of trace class operators, one can 
still define the Chern character via the secondary characteristic classes. The essential reason for this is that along vertical directions in $\pi: \Cal A \to \Cal A/\Cal G$
the connection $\tilde\theta$ agrees with $\omega(\cdot ; \theta)$ and the latter takes values in the Lie algebra of $U^1.$ For example, a spherical cycle
$S^{2k} \subset \Cal A/\Cal G$ can be lifted to $\Cal A$ after removing one point $x$ from $S^{2k}.$ The inverse image $\pi^{-1}(x) \subset \Cal A$ is then a cycle
of dimension $2k-1$ inside of a gauge orbit. The integral of the ill-defined form $\text{tr} \, \Cal F^k$ over the cycle $S^{2k}$ becomes then an integral
of the corresponding Chern-Simons form over the cycle in the fiber $\pi^{-1}(x).$ But this integration involves only the vertical directions in the connection
and curvature, and are thus well-defined because of the trace class property of the connection along vertical directions. This process can be generalized 
to arbitrary homology cycles as explained below.

Recall that in the trace class case, on contractible sets $c^{2n} := \tr\, \Cal F^n= d \tr\, C^{2n-1}
= d c^{2n-1}$
where $C^{2n-1}$ is a differential polynomial of the connection form.  Explicitly, in a general vector bundle with a local connection 1-form $A$ and the corresponding
curvature form $F$ the Chern-Simons form $c^{2n-1}$ is given by the formula
$$  c^{2n-1} = n \int_0^1 dt \, \text{tr}\, A_t F_t^{n-1}\tag2.9$$
where $A_t = tA$ and $F_t = t dA + \frac{t^2}{2} [A,A].$

The integral of $c^{2n}$ over a simplicial cycle $c_{2n}$ of dimension $2n$ is then
given as the sum of integrals of $c^{2n-1}$ over all the faces $c_{2n-1, \alpha\beta}$ in the singular simplex $c_{2n}= \sum c_{2n, \alpha}.$ The convention
is that on the face $c_{2n+1, \alpha\beta}$ one integrates the pull-back form $c^{2n-1}_{\alpha}$ (defined by the local section $\psi_{\alpha}$) over the boundary
component $c_{2n-1, \alpha\beta}= c_{2n,\alpha} \cap c_{2n,\beta} $ with the orientation defined by the outward normal to $c_{2n,\alpha}$ and the fixed orientation
on $c_{2n}.$ For a given index pair $\alpha,\beta$ one has the pair of integrals
$$ \int_{c_{2n-1, \alpha\beta}}  c^{2n-1, \alpha}    \text{  and } \int_{c_{2n-1,\beta\alpha}} c^{2n-1, \beta}.$$  
Now the form $c_{2n-1, \alpha\beta}$ is a gauge transform of $c_{2n-1, \beta\alpha}$ with respect to the change of local trivialization $\psi_{\beta} = \psi_{\alpha} \cdot 
g_{\alpha\beta}$ so, taking into account the orientations, the above pair gives the contribution
$$ \int_{c_{2n-1, \alpha\beta}} ( c^{2n-1, \alpha}  -  g^*_{\alpha\beta} c^{2n-1, \alpha} )  = \int_{c_{2n-1,\alpha\beta}} c^{2n-1, \alpha\beta}.$$  
The difference under the integral is a trace of a differential polynomial in the local connection form $\psi_{\alpha}^* \tilde \theta$ and
 the Maurer-Cartan 1-form
$g_{\alpha\beta}^{-1} dg_{\alpha\beta};$ using (2.9) one observes that each term at least of degree one in the latter. When the Maurer-Cartan form takes values in the Lie algebra
of $U^1,$ that is, in the ideal of trace-class operators, all the traces actually converge.  \enddemo

Thus we can conclude, using Corollary 1  and  Lemma 3,

\proclaim{Theorem 2}  Let $\Cal F$ be the curvature corresponding to the connection $\nabla,$ defined from the 1-form $\tilde\theta$ on $\Cal A.$
The formal Chern character $\text{str}_c\, e^{\Cal F/2\pi i} $   defines a cohomology class on the base $X$ as a pairing with singular cycles, evaluated 
using the descent equations in the {\v C}ech - de Rham double complex up to the cocycle in $H^{(*, 1)}$ of {\v C}ech degree one.  The class does not depend on the choice of the regularization $T_A.$    \endproclaim

\demo{Proof} What remains to show is that the Chern character in cohomology does not depend on the choice of the family $T_A.$ A different choice $T'_A$ is 
related to the first by $T'_A = T_A S_A$ with $S_A \in U_1.$ Since the space $\Cal A$ of smooth connections is flat, we can choose a homotopy $T_A(t) =
T_A S_A(t)$ with $S_A(0) = 1$ and $S_A(1) = S_A.$ This gives a homotopy $\nabla(t)$ of the connections, connecting $\nabla$ to the connection
$\nabla'$ defined by $T_A'.$ The rest is then the standard argument relating the Chern characters defined by a pair of connections in a vector bundle, the
difference being the exterior derivative of a Chern-Simons form. 

One should also note that the choice of the base point $\epsilon$ is not essential. Any other base point $\epsilon'$ (a self-adjoint operator with ${\epsilon'} ^2 =1$ and
essential spectrum at $\pm 1$) is written as $\epsilon' = W \epsilon W^{-1}$ for some unitary operator $W.$ The connection $W\tilde\theta W^{-1}$ satisfies then the
same conditional 
trace-class properties (in the vertical directions) with respect to $\epsilon'$ as does $\tilde\theta$ with respect to $\epsilon$ and so the Chern character defined by the pair $(\epsilon, \tilde\theta)$ is equal to the Chern character defined by
$(\epsilon', \tilde\theta').$   \enddemo $\square$

\vskip 0.3in
3. THE ODD CASE: CHERN CHARACTER FOR $U_p$ BUNDLES

\vskip 0.2in 

When $M$ is an odd dimensional spin manifold and $D$ a family of Dirac operators acting on sections of the tensor product of the spin bundle
and a (trivial) vector bundle over $M,$ the Chern character takes values in the odd cohomology groups of the parameter space $X.$ 
In the following I want to describe how the renormalization method in Section 2,  with the operators like $T_A,$ can be used to compute
the Chern character on $X=\Cal A/ \Cal G.$

The first step is to replace the family of Dirac operators $D_A$ by a family of operators $\epsilon + B$ where $B$ is a bounded operator (which depends on
$A$) such that $[\epsilon, B] \in L_2$ and $\epsilon B + B\epsilon \in L_1.$ First, choose a  family of operators $T_A$  such that $T_A^{-1} \epsilon T_A -F_A$ is
trace-class, [LM]; equivalently, $\epsilon - {T_A}^{-1} F_A T_A$ is trace-class.  Actually, since the operators $D_A$ might have zero modes we have to replace the
sign operator $F_A$ by an approximate sign operator, say $F_A = D_A (|D_A| + \exp(-D_A^2))^{-1}.$ Since these operators transform equivariantly with respect to gauge 
transformations, they descend to a family of bounded operators parametrized by $X= \Cal A/\Cal G.$  We denote $T_A F_A T_A^{-1} = \epsilon + B(A).$

As in the even case, we have a cocycle $\Omega: \Cal A \times \Cal G \to U_1.$ Since in the odd case we do not have the grading $\Gamma,$ there is no
further refinement of this cocycle to a $U^1$ valued cocycle. The cocycle  $\omega$ defines a $U_1$ bundle $P$ over the moduli space $X.$  The Atiyah-Singer
connection in the principal bundle $\Cal A \to X$ defines again a connection $\nabla$ in the $U_1$ bundle (and in the associated Hilbert bundle) over $X.$ 
We  can then define a \it superconnection \rm 

$$\Theta = \delta + B + \nabla.$$

Here $\delta$ is an exterior differentiation in the sense of noncommutative geometry, [Co]: It is acting on 'k-forms' of the type $\Phi= \alpha_0 \delta \alpha_1 \dots \delta \alpha_k$
as $\delta\Phi = \delta \alpha_0\delta\alpha_1 \dots \delta \alpha_k$ where the $\alpha_i$'s are bounded operators in $H$ such that $\delta\alpha = [\epsilon, \alpha] \in
L_p.$  One can check that when $k$ is odd $\delta\Phi = \epsilon\Phi + \Phi\epsilon$ whereas when $k$ is even $\delta \Phi = [\epsilon, \Phi].$  In this sense, $B$ is a 
1-form.  We have a double complex of forms $\Omega^{k,l}$ where $k$ gives the $\delta$-degree and $l$ the form degree related to the exterior differentiation on the
base $X.$ Finally, we modify the the action $\delta \mapsto (-1)^l \delta$ in order that $\delta d + d \delta =0.$ 

Choosing an open cover $\{U_{\alpha}\}$ of $X$  with local trivializations of the bundle $P$ and transition functions $\tilde g_{\alpha\beta}:U_{\alpha\beta} \to U_1,$
the local operator families $B_{\alpha}$ are related by
$$B_{\beta} = \tilde g_{\alpha\beta}^{-1} B_{\alpha} \tilde g_{\alpha\beta} + \tilde g_{\alpha\beta}^{-1}[\epsilon, \tilde g_{\alpha\beta}]$$
and the local 1-forms $\eta_{\alpha}$ on $X$ which define the connection $\nabla$ locally as $\nabla = d + \tilde\theta_{\alpha}$ are related by
$$ \tilde\theta_{\beta} = \tilde g_{\alpha\beta}^{-1} \tilde\theta_{\alpha} \tilde g_{\alpha\beta} + \tilde g_{\alpha\beta}^{-1} d \tilde g_{\alpha\beta}$$
In this way, the superconnection transforms covariantly,
$$\Theta_{\beta} = \tilde g_{\alpha\beta}^{-1}\Theta_{\alpha} \tilde g_{\alpha\beta}.$$

When the $U_1$ bundle arises as an associated bundle to the gauge bundle $\Cal A \to \Cal A/\Cal G$ through the cocycle $\Omega$ a superconnection
can be defined with the help of local operator families $\epsilon + B_{\alpha}(A)  = F_A$ defined using a local
section $\psi_{\alpha}:  \Cal A/\Cal G\to \Cal A.$ 

The square of the superconnection gives the graded curvature,
$$\Theta^2  =   \Cal F^{(2,0)} + \Cal F^{(1,1)} + \Cal F^{(0,2)}$$
with $\Cal F^{(2,0)} = \epsilon B + B\epsilon + B^2, \, \Cal F^{(1,1)} = [\epsilon, \tilde\theta] +  d B + [B,\tilde\theta]$ and $\Cal F^{(0,2)} =  d \tilde\theta + \tilde\theta^2.$

Since the classifying space  $BU_1=U(\infty)$ has odd homotopy type, the characteristic classes of $U_1$ bundles are located in
in odd cohomology of the base space. The lowest nontrivial form is in $H^1(X).$  In a similar way as in the even case, we have:

\proclaim{Theorem 3} The terms $\tr_c \ (\Cal F^n)^{(1,2n-1)}$  in the formal Chern character on $X$ define odd closed  index forms on $X= \Cal A/\Cal G,$ 
through the terms of {\v C}ech degree one in the corresponding {\v C}ech - de Rham complex.   
They do not depend on the choice of the
regularization family $T_A.$ \endproclaim

As in the even case,  if one chooses  a connection form $\tilde\theta$ such that $[\epsilon, \tilde\theta]$ is trace-class then also the forms $\text{tr} \, \Cal F^n$ are 
well-defined directly. 

Let us consider two simple examples, where we can in fact choose $T_A \equiv 1.$ 

\bf Example 1 \rm For  $X=S^1$ the nontrivial $U_1$ bundles are
classified by the homotopy type of the transition function $g\in U_1$ at the ends of the interval $0\leq t \leq 2\pi$ parametrizing the circle
$e^{ i t}\in S^1.$ The group consists of disconnected components labelled by the Fredholm index of the $a$ block in
$$g=\left(\matrix a& b\\ c& d\endmatrix\right).$$

The  superconnection is given by $\Theta= \delta  + d + \frac{t}{2\pi} g^{-1}[\epsilon,g]$  and $\tilde\theta=0.$  
The representative for the index form in $H^1(S^1)$ is then $\frac{1}{2\pi} dt\, \tr \,g^{-1}[\epsilon,g]=  \tr\, \Cal F^{(1,1)}.$ The integral of this form over the circle $S^1$
is equal to twice the Fredholm index of $a.$

\bf  Remark \rm Actually, it is sufficient that the off-diagonal blocks of all components $\Cal F^{(a,b)}$ are Hilbert-Schmidt (so we can take as the structure group of the bundle $U_2$ instead of
$U_1$) and in addition the diagonal blocks of $\Cal F^{(1,1)}$ are trace-class.
Then all the odd forms on $X$ arising in $\tr\, \Cal F^k$ are conditionally converging traces.

\bf Example 2 \rm  The space $\Cal A$ of smooth vector potentials on the circle $S^1$ with values in the Lie algebra $\frak g$ of the gauge group
$G= SU(2)$ is a principal $\Omega G$ bundle over $G.$ The group $\Omega G$ of based smooth loops acts freely, as gauge transformations, on
$\Cal A.$  The group $\Omega G$ acts in the Hilbert space $H$ of square integrable functions on the circle with values in $\Bbb C^2$ 
by point-wise multiplication. The space $H= H_+\oplus H_-$ is polarized to nonnegative and negative Fourier modes. It is known
that $\Omega G$ embeds to $U_2(H),$ [PS]. (In fact, also $\Omega G \subset U_1$).
 Thus the $\Omega G$ bundle over $SU(2)$ extends naturally to a $U_2$ bundle. 
A supeconnection is defined as follows.  On the equator $S^2 \subset S^3 = SU(2)$ we have the transition function $g: S^2 \to \Omega G$
(we extend the equator to a small open neighborhood thereof and extend $g$ as a constant to the normal direction). 
Explicitly, $g$ can be given as 
$$ g(\bold n) =  \exp(i f(\phi) \bold n \cdot \bold \sigma) $$
where $\bold n = (n_1, n_2, n_3)$ is a unit vector in $\Bbb R^3$ parametrizing the sphere $S^2$ and $\bold \sigma = (\sigma_1, \sigma_2, \sigma_3)$ is the set 
of complex hermitean trace-less  $2\times 2$ Pauli matrices, $\sigma_k^2= 1$ and $\sigma_1 \sigma_2 = i \sigma_3,$ and cyclic permutations of this relation.
$f(\phi)$ is a smooth function on the interval $[0,2\pi]$ with $f(0) =0$ and $f(2\pi) = n \cdot 2\pi,$  $n\in \Bbb Z.$
We set $B=0= \tilde\theta$ on the half sphere $x_3 > 0$ and $B= \frac{t}{\pi} [\epsilon, g]g^{-1},
\tilde\theta= \frac{t}{\pi} g^{-1}dg$ on the other hemisphere, where $0\leq t \leq \pi$ is the distance from the 'South Pole' on $S^3.$  This defines the superconnection
$\Theta$ on the Southern hemisphere whereas $\Theta = \delta + d$ on the Northern hemisphere.  The $(1,3)$ component of the square of
the supercurvature on the Southern hemisphere is now
$$\tr_H \, (\Cal F^2)^{(1,3)} = \frac{dt}{\pi} ( - \frac{t}{\pi} +(\frac{t}{\pi})^2)   \tr_H \, g^{-1}[\epsilon,g] (g^{-1}dg)^2
$$
and the integral of the square of the curvature over the whole  3-sphere becomes
$$ I(g)= \frac{i}{12\pi } \int_{0\leq \phi \leq 2\pi} \int_{S^2} \tr_{\Bbb C^2}  \, (g' g^{-1}) (dg g^{-1})^2 $$
where $g'$ is the derivative with respect to $\phi$ and $d$ stands for the derivatives to the $S^2$ directions. We have used the formula
$$ \tr_H  [\epsilon, X] Y = \frac{1}{2\pi i}  \int_{S^1} \tr_{\Bbb C^2}  \, X' Y d\phi$$
for matrix valued functions on the unit circle. The value of the integral $I(g)$ is equal to $2\pi i$ times the winding number of the map $\tilde g: S^3 \to G,$
where $\tilde g$ is defined as $\tilde g(\phi, \bold n) = g(\bold n)(\phi)$ with $\bold n\in S^2$ and $\phi$ parametrizes $S^1;$ since $g$ is a map from $S^2$ to the \it based
loop group \rm $\Omega G,$ this function on $S^2 \times S^1$ lifts to a map on $S^3.$ For the particular choice $g$ above, the winding number is $2n.$ 
With a minor modification one can get  also the odd winding numbers: For that purpose one defines $g(\bold n)(\phi)$ for $0\leq \phi \leq \pi$ as above, but
for $\pi \leq \phi \leq  2\pi$ one sets $g(\bold n)(\phi) = h(\phi)$ where $h$ is a fixed path in $SU(2)$ with end points $\pm \bold 1.$ Now the winding number of $\tilde g$
is equal to $n.$
\vskip 0.3in
REFERENCES
\vskip 0.3in

[AS1] M.F. Atiyah and I. M. Singer: The index of elliptic operators. IV. Ann. of Math. (2) \bf 93 \rm (1971), 119 - 138.

[AS2] M.F. Atiyah and I.M. Singer:  Dirac operators coupled to vector potentials. Proc. Nat. Acad. Sci. U.S.A. \bf 81 \rm (1984), no. 8, Phys. Sci., 2597 - 2600.

[B] J.-M.  Bismut: The Atiyah-Singer index theorem for families of Dirac operators: two heat equation proofs. Invent. Math. \bf 83 \rm (1985), no. 1, 91 - 151.

[CMM]  A. Carey, J. Mickelsson, and M. Murray:  Index theory, gerbes, and Hamiltonian quantization. Comm. Math. Phys. \bf 183 \rm (1997), no. 3, 707 - 722.
 
[C] Alain Connes: \it Noncommutative Geometry. \rm Academic Press (1994); also http://www.alainconnes.org/docs/book94bigpdf.pdf

[K]  N.  Kuiper:  The homotopy type of the unitary group of Hilbert space. Topology \bf 3 \rm (1965),  19 - 30.

[LM] E. Langmann and J. Mickelsson: Scattering matrix in external field problems. J. Math. Phys. \bf 37 \rm(1996), no. 8, 3933 - 3953. 

[M] J. Mickelsson:  Wodzicki residue and anomalies of current algebras. In: Integrable models and strings (Espoo, 1993), 
123 - 135, Lecture Notes in Phys., \bf 436, \rm Springer, Berlin, (1994).

[Mi] John Milnor:  On infinite dimensional Lie groups. Preprint, Institute for Advanced Studies, Princeton (1982)

[MR] J. Mickelsson and S. Rajeev:  Current algebras in d+1-dimensions and determinant bundles over infinite-dimensional Grassmannians. 
Comm. Math. Phys. \bf 116 \rm  (1988), no. 3, 365 - 400.

[P] Richard Palais:  On the homotopy type of certain groups of operators. Topology \bf 3 \rm (1965),  271 - 279.

[PS] A. Pressley and G. Segal: \it Loop Groups. \rm  Clarendon Press, Oxford (1986)

[Q] Daniel Quillen: Superconnections and the Chern character. Topology \bf 24 \rm (1985), no. 1, 89 - 95.

\enddocument